\documentclass[letterpaper, 10 pt, conference]{ieeeconf}

\newcommand{\myreferences}{compbib.bib}

\IEEEoverridecommandlockouts                              
\usepackage{enumerate}
\usepackage{amsmath,amssymb,algorithm,cite,cases,bm,float,subfigure,graphicx,url}
\usepackage{pbox}
\newtheorem{defn}{Definition}[section]

\newtheorem{thm}{Theorem}[section]
\newtheorem{lem}{Lemma}[section]

\usepackage{amsmath,amssymb,algorithm,cite,subfigure,caption,cases,bm,float,subfigure,graphicx,url,color}
\usepackage{comment}
\usepackage{thm-restate, enumerate }

\definecolor{darkred}{RGB}{150,0,0}
\definecolor{darkgreen}{RGB}{0,150,0}
\definecolor{darkblue}{RGB}{0,0,200}

\newcommand{\labe}{(\la_\beta)}

\newcommand{\wt}{\tilde{\mathbf{w}}_{\sigma}}

\newcommand{\dd}{\tilde d}

\newcommand{\rP}{\xrightarrow{P}}

\newcommand{\Rc}{\mathcal{R}}
\newcommand{\R}{\mathbb{R}}

\newcommand{\lims}{\lim_{\sigma\rightarrow 0}}

\newcommand{\phiz}{\phi_0}

\newcommand{\taubb}{\tau_{\text{best}}}

\newcommand{\cost}{\mathrm{Cost}}

\newcommand{\nlim}{\lim_{n\rightarrow\infty}}

\newcommand{\NSE}{\operatorname{NSE}}
\newcommand{\aNSE}{\operatorname{aNSE}}
\newcommand{\wNSE}{\operatorname{wNSE}}

\newcommand{\La}{\Lambda}

\newcommand{\D}{{D(\tau)}}
\newcommand{\C}{{C(\tau)}}

\newcommand{\beq}{\begin{equation}}
\newcommand{\eeq}{\end{equation}}
\newcommand{\bea}{\begin{align}}
\newcommand{\eea}{{\end{align}}}

\newcommand{\M}{\mathbf{M}}

\newcommand{\labb}{\lambda_{\text{best}}}

\newcommand{\G}{\mathbf{G}}

\newcommand{\A}{\mathbf{A}}

\newcommand{\bu}{\boldsymbol{\pi}}
\newcommand{\buh}{\boldsymbol{\pi}_{\la\paf}(\h)}

\newcommand{\dtth}{\mathrm{dist}^2(\h,{\la\paf})}

\newcommand{\w}{{\mathbf{w}_\sigma}}
\newcommand{\ww}{{\mathbf{w}}}

\newcommand{\x}{\mathbf{x}}
\newcommand{\ub}{\mathbf{u}}
\newcommand{\g}{\mathbf{g}}
\newcommand{\vb}{\mathbf{v}}
\newcommand{\bb}{\mathbf{b}}

\newcommand{\y}{\mathbf{y}}
\newcommand{\s}{\mathbf{s}}
\newcommand{\z}{\mathbf{z}}

\newcommand{\ab}{\mathbf{a}}

\newcommand{\h}{\mathbf{h}}

\newcommand{\wh}{\hat{\mathbf{w}}_\sigma}

\newcommand{\Sc}{{\mathcal{S}}}
\newcommand{\map}{{\mathrm{map}}}

\newcommand{\Lc}{{\mathcal{L}}}

\newcommand{\Nn}{\mathcal{N}}
\newcommand{\Cc}{\mathcal{C}}

\newcommand{\Pro}{{\mathbb{P}}}
\newcommand{\E}{{\mathbb{E}}}

\newcommand{\paf}{\pa f(\x_0)}

\newcommand{\la}{{\lambda}}

\newcommand{\eps}{\epsilon}

\newcommand{\pa}{\partial}

\newcommand{\dt}{\text{{dist}}}

\newcommand{\vs}{\vspace}

\newcommand{\nn}{\nonumber}

\begin{document}
\title{\LARGE{\bf{Asymptotically Exact Error Analysis for the Generalized $\ell_2^2$-LASSO  }}\vs{1pt}}
%
\author{
Christos Thrampoulidis$^*$, Ashkan Panahi$^\dagger$ and Babak Hassibi$^*$\vs{2pt}\\
Department of Electrical Engineering, Caltech, Pasadena -- 91125
\thanks{
$^*$Department of Electrical Engineering, Caltech, Pasadena \newline
$^\dagger$Signal Processing Group, Chalmers Univ. of Technology, Gothenburg
}
}

%

\maketitle
\begin{abstract} 
Given an unknown signal $\mathbf{x}_0\in\mathbb{R}^n$ and linear noisy measurements $\mathbf{y}=\mathbf{A}\mathbf{x}_0+\sigma\mathbf{v}\in\mathbb{R}^m$, the generalized $\ell_2^2$-LASSO solves $\hat{\mathbf{x}}:=\arg\min_{\mathbf{x}}\frac{1}{2}\|\mathbf{y}-\mathbf{A}\mathbf{x}\|_2^2 + \sigma\lambda f(\mathbf{x})$. Here,  $f$ is a convex regularization function (e.g. $\ell_1$-norm, nuclear-norm) aiming to promote the structure of $\mathbf{x}_0$ (e.g. sparse, low-rank), and, $\lambda\geq 0$ is the regularizer parameter. A related optimization problem, though not as popular or well-known, is often referred to as the generalized $\ell_2$-LASSO and takes the form $\hat{\mathbf{x}}:=\arg\min_{\mathbf{x}}\|\mathbf{y}-\mathbf{A}\mathbf{x}\|_2 + \lambda f(\mathbf{x})$, and has been analyzed in \cite{OTH}. \cite{OTH} further made conjectures about the performance of the generalized $\ell_2^2$-LASSO. This paper establishes these conjectures rigorously. We measure performance with the normalized squared error $\mathrm{NSE}(\sigma):=\|\hat{\mathbf{x}}-\mathbf{x}_0\|_2^2/\sigma^2$. Assuming the entries of $\mathbf{A}$ and $\mathbf{v}$ be i.i.d. standard normal,  we precisely characterize the ``asymptotic NSE" $\mathrm{aNSE}:=\lim_{\sigma\rightarrow 0}\mathrm{NSE}(\sigma)$ when the problem dimensions $m,n$ tend to infinity in a proportional manner. The role of $\lambda,f$ and $\mathbf{x}_0$ is explicitly captured in the derived expression via means of a single geometric quantity, the Gaussian distance to the subdifferential. We conjecture that $\mathrm{aNSE} = \sup_{\sigma>0}\mathrm{NSE}(\sigma)$. We include detailed discussions on the interpretation of our result, make connections to relevant literature and perform computational experiments that validate our theoretical findings.


\end{abstract}

\section{Introduction}

\subsection{Generalized LASSO}
The Generalized $\ell_2^2$-LASSO has  emerged as a powerful tool for the recovery of \emph{structured} signals (sparse, low rank, etc.) from linear noisy measurements in a variety of applications in statistics, signal processing, machine learning, etc.. Given an unknown signal $\x_0\in\R^n$ and measurements $\y=\A\x_0+\z\in\R^m$, it solves:
 \begin{align}\label{eq:intro_GL}
 \hat\x:=\arg\min_{\x}({1}/{2})\|\y-\A\x\|_2^2 + \sigma\la f(\x).
 \end{align}
Here,  $f$ is a \emph{convex} regularization function, typically non-smooth (e.g. $\ell_1$-norm, nuclear-norm, $\ell_1/\ell_2$-norm), aiming to promote the structure of $\x_0$ (e.g. sparse, low-rank, block-sparse). $\la\geq 0$ is the regularizer parameter and is scaled with the standard deviation $\sigma$ of the noise vector $\z$, which is typically modeled to have entries i.i.d. $\Nn(0,\sigma^2)$. The term ``LASSO" was coined by Tibshirani \cite{TibLASSO} who first introduced \eqref{eq:intro_GL} with $f$ chosen as the $\ell_1$-norm. In this view, \eqref{eq:intro_GL} is a natural generalization to other structures and \emph{convex} regularizers. We have added the indicator ``$\ell_2^2$" to distinguish \eqref{eq:intro_GL} from a variant which takes the form \cite{Belloni,OTH}:
\vspace{-5pt}
 \begin{align}\label{eq:intro_GLOTH}
 \hat\x:=\arg\min_{\x}\|\y-\A\x\|_2 + \mu f(\x).
 \end{align}
 We call this the Generalized $\ell_2$-LASSO, but it is also known in related literature (e.g. \cite{Belloni}) as the square-root LASSO. The two optimizations in \eqref{eq:intro_GL} and  \eqref{eq:intro_GLOTH} are fundamentally related: from  optimality conditions there exists a \emph{mapping} between the regularizer parameters $\la$ and $\mu$ for which the performance is equivalent. However, not only is this mapping non-trivial to characterize, but also there exist other differentiating features. For instance, note that in \eqref{eq:intro_GLOTH} the regularizer parameter $\mu$ need not scale (thus is agnostic) with the noise variance \cite{Belloni, OTH}. A comparison between the two algorithms is beyond the scope of the paper, but our result, when combined with those of \cite{OTH}, inevitably results in some further related discussions in the next sections.  In what follows, we often drop the attribute ``Generalized" and simply refer to \eqref{eq:intro_GL} and \eqref{eq:intro_GLOTH} as the $\ell_2^2$-LASSO and $\ell_2$-LASSO, respectively.
 
 \subsection{Performance Analysis and Related Literature}
 
 A natural measure of performance of \eqref{eq:intro_GL} or \eqref{eq:intro_GLOTH} is the Normalized Squared Error
 $
 \NSE:= {\|\hat\x - \x_0\|_2^2}/{\sigma^2}.\nn
$
To facilitate the theoretical analysis of the $\NSE$, it is standard to assume that the measurement matrix $\A$ is drawn at random from some ensemble. Early well-known bounds on the $\NSE$  were  order-wise in nature (i.e. accurate only up to constant multiplicative factors) and derived based on RIP and Restricted Eigenvalue assumptions on the measurement matrix \cite{candes2006stable,candes2007dantzig,bickel,negahban2012unified,Belloni}. To the best of our knowledge, the first precise formulae predicting the limiting behavior of the $\ell_2^2$-LASSO reconstruction error were provided by Donoho, Maleki, and Montanari \cite{arianSensitivity}; a proof appeared lated by Bayati and Montanari in \cite{MontanariLASSO}. The authors of these references consider the $\ell_2^2$-LASSO with $\ell_1$-regularization, i.i.d Gaussian sensing matrix $\A$ and use the Approximate Message Passing (AMP) framework for the analysis (also see subsequent related works \cite{arianComplex,arianDenoising}).
 More recently, Stojnic \cite{StoLASSO} introduced an alternative framework and used it to derive a tight upper bound on the NSE of the following constrained version of the LASSO:
 \vspace{-5pt}
 \begin{align}\label{eq:intro_CL}
 \min_\x{\|\y-\A\x\|_2}\quad\text{s.t.}\quad \|\x\|_1\leq \|\x_0\|.
 \end{align}
 Stojnic's approach cleverly uses a comparison lemma due to Gordon \cite{gorLem} , known as the Gaussian min-max Theorem (GMT). What allowed him to use this machinery in the first place was the observation
\footnote{In fact, this same trick is used in the classical application of GMT that lower bounds the minimum singular value of an i.i.d. Gaussian matrix, \cite{vershynin2010introduction}.}
  that \eqref{eq:intro_CL} can be equivalently expressed as a min-max problem as follows:
\vspace{-5pt}
 \begin{align}\label{eq:intro_Sto}
 \min_\x\max_{\|\ub\|\leq 1}{\ub^T(\y-\A\x)}\quad\text{s.t.}\quad \|\x\|_1\leq \|\x_0\|_1.
 \end{align} 
 It turns out that this form is appropriate for the application of GMT. The same idea was used in \cite{OTH} to generalize the results of \cite{StoLASSO} to arbitrary convex regularizer functions in \eqref{eq:intro_CL}. However, the main contribution of \cite{OTH} is the extension of the results to the generalized $\ell_2$-LASSO. The presence of the regularizer parameter $\la$ in \eqref{eq:intro_GLOTH} makes the extension 
   non-trivial and considerable effort had to be undertaken in \cite{OTH}. Of course, the same observation that allows the use of the GMT in the first place, is here the same as in \eqref{eq:intro_Sto}, namely \eqref{eq:intro_GLOTH} can be expressed as
   \vspace{-5pt}
 \begin{align}\label{eq:intro_Sto}
 \min_\x\max_{\|\ub\|\leq 1}{\ub^T(\y-\A\x)}+\mu f(\x).
 \end{align}  
At that time it wasn't clear to the authors of \cite{OTH} how to leverage the objective function in \eqref{eq:intro_GL} and analyze the $\NSE$ of the $\ell_2^2$ LASSO under the same machinery.  However, making an ``educated guess" on the formula that governs the mapping between the two versions of the LASSO, they were able to translate results from \eqref{eq:intro_GLOTH} to \eqref{eq:intro_GL}. This led them to \emph{conjecture} a formula for the upper bound on the $\NSE$ of the $\ell_2^2$-LASSO, which was also suggested by numerical simulations. 

\vspace{-2pt}
\subsection{Our Contribution}
In this work, we \emph{rigorously} establish the conjecture raised in \cite{OTH} on the $\NSE$ of the Generalized $\ell_2^2$-LASSO under i.i.d. Gaussian measurements. Instead of worrying about the mapping function between \eqref{eq:intro_GL} and \eqref{eq:intro_GLOTH} and translating the results from the latter to the former, we follow a direct approach. The \emph{key observation} is that the objective function in \eqref{eq:intro_GL} can be appropriately linearized for the purpose of using the GMT, and be written equivalently as:
\vspace{-5pt}
 \begin{align}\nn
 \min_\x\max_{\ub}\ub^T(\y-\A\x)-({1}/{2})\|\ub\|^2+\la\sigma f(\x).
 \end{align}  
Beyond this trick, what facilitates our analysis is a result from \cite{tight}. Essentially, \cite{tight} builds a clear, concrete and easy to apply framework based on Stojnic's original idea of combining GMT with convexity. This allows a more insightful and compact analysis when compared to \cite{StoLASSO,OTH}.


\section{Result}

\subsection{ Setup}

Let $\x_0\in R^n$, $\y= \A\x_0 + \sigma\vb\in \R^m$ and \emph{convex} $f:\R^n\rightarrow\R$. The $\ell_2^2$-LASSO solves \eqref{eq:intro_GL}
for $\la\geq 0$. The reconstruction vector $\hat\x$ depends explicitly on $\A,\la,\sigma, f$, and, implicitly on $\vb,\x_0$ through the measurement vector $\y$.
Define the Normalized Squared-Error of \eqref{eq:intro_GL} as 
\begin{align}\label{eq:NSE_def}
\NSE(\sigma) := {\|\hat\x-\x_0\|_2^2}/{\sigma^2}.
\end{align}

\subsubsection{Assumptions}\label{sec:ass}
We assume that the entries of $\A$ and $\vb$ are i.i.d. $\Nn(0,1)$. The regularizer $f:\R^n\rightarrow\R$ is  convex and continuous.
Also, $\x_0$ is not a minimizer of $f$. Popular regularizers include the $\ell_1$-norm, nuclear-norm, $\ell_{1,2}-norm$ etc. (please refer to \cite{Cha,TroppEdge} for further examples).

\subsubsection{Large system limit}\label{sec:large}
Our results hold in an asymptotic regime in which the problem dimensions grow to infinity. We consider  a sequence of problem instances $\{\A,\vb,\x_0,f\}_{m,n}$ as in \eqref{eq:intro_GL} indexed by $m$ and $n$ such that both $m,n\rightarrow\infty$. In each problem instance, $\A, \vb$ and $f$ satisfy the assumptions of Section \ref{sec:ass}. Furthermore, $\hat\x$ and $\NSE(\sigma)$ denote the output of \eqref{eq:intro_GL} and the corresponding NSE. To keep notation simple, we avoid introducing explicitly the dependence of variables on the problem dimensions $m,n$. 

\subsubsection{NSE: worst-case and asymptotic}
Define the \emph{worst-case} NSE as
$
\wNSE := \sup_{\sigma>0} \NSE(\sigma).
$ We say that recovery of $\x_0$ by means of \eqref{eq:intro_GL} is \emph{robust} whenever $\wNSE<\infty$. Further, define 
 the \emph{asymptotic} NSE as
$
\aNSE := \lim_{\sigma\rightarrow 0} \NSE(\sigma).
$
Theorem \ref{thm:main} in Section \ref{sec:main} derives a precise expression for $\aNSE$ in the large system limit. In  Section \ref{sec:rem} we conjecture that under our assumptions $\aNSE=\wNSE$, which highligths  the significance of studying the $\aNSE$. Recent results \cite{Don94,arianSensitivity,oymakProx,StoLASSO,OTH}, have shown that $\wNSE$ is also achieved in the limit $\sigma^2\rightarrow 0$ for algorithms of nature similar to \eqref{eq:intro_GL} under similar setups. Please also refer to relevant discussion (on the similarly defined notion of \emph{noise-sensitivity}) in \cite{Yu}.

\subsubsection{Gaussian Squared Distance}   \mbox{The subdifferential of $f$ at $\x_0$ is the set of vectors:} \mbox{
$\paf = \left\{ \s\in\mathbb{R}^n | f(\x_0+\ub) \geq f(\x_0) + \s^T\ub, \forall \ub\in\mathbb{R}^n  \right\}.$}
It is nonempty, convex and compact \cite{Roc70}. Also, it does not contain the origin (recall $\x_0$ is not a minimizer).
For any nonnegative number $\tau\geq0$, denote the scaled (by $\tau$) subdifferential set as $\tau\paf=\{\tau\s|\s\in\paf\}$. Also, for the conic hull of the subdifferential $\paf$, write $\text{cone}(\paf) = \{\s | \s\in\tau\paf, \text{ for some } \tau\geq 0\}$. 
For $\Cc\subset\mathbb{R}^n$ nonempty, convex, closed set and $\ub\in \R^n$, denote the projection and distance 
as $\bu_\Cc(\h):=\arg\min_{\s\in\Cc}\|\ub-\s\|_2$ and
$\dt(\ub,\Cc) = \|\ub-\bu_\Cc(\h)\|_2$.

\begin{defn}[Gaussian squared distance] Assume $f:\R^n\rightarrow\R$ convex. Let $\h\in\mathbb{R}^n$ have i.i.d $\Nn(0,1)$ entries. The gaussian squared distance to the scaled subdifferential is defined as
\begin{align}\label{eq:D_def}
\D:= D_{\paf}(\tau) := \E_\h\left[\dtth\right].
\end{align}
\end{defn}
$\D$ appears as a fundamental quantity in the study of the phase transitions of noiseless compressive sensing: it has been shown  that 
\begin{align}\label{eq:m>D}
m\gtrsim\min_{\tau\geq 0}D(\tau)\approx\E[\dt^2(\h,\text{cone}(\paf)].
\end{align}
 is sufficient \cite{Sto,Cha} and necessary \cite{TroppEdge} for the recovery of $\x_0$ from \emph{noiseless} linear observations. Thus, it is no surprise that the properties of $\D$ have been analyzed in detail in\cite[Lem.~C.2]{TroppEdge} (also, \cite[Lem.~8.1]{OTH}). The same quantity plays central role in the analysis of the noisy case considered here; we make precise reference to relevant properties whenever they appear useful throughout our exposition. For the statement of our results, we need the following: $\D$ is differentiable for $\tau>0$ and $\partial D(\tau)/\partial\tau = -(2/\tau) \C$,
\begin{align}\label{eq:C_def}
\C := \E_\h\left[(\h-\buh)^T\buh\right].
\end{align}

To familiarize with the definitions in \eqref{eq:D_def} and \eqref{eq:C_def}, it is instructive to specialize to the case where $f$ is the $\ell_1$-norm and $\x_0$ is a $k$-sparse vector. Then, $\paf$ has a simple characterization and $\D,\C$ admit simple closed-form expressions in terms of the tail distribution $Q(\tau)$ of a standard Gaussian (e.g.,\cite[App.~H]{OTH} ):
\begin{align}
D(\tau) &= k(1+\tau)^2+(n-k)(2(1+\tau^2)Q(\tau)-\sqrt{\frac{2}{\pi}}\tau e^{-\frac{\la^2}{2}}) \nn\\
C(\tau) &= -k\tau^2+(n-k)(2\tau^2 Q(\tau)-\sqrt\frac{2}{\pi}\tau e^{-\frac{\tau^2}{2}})\label{eq:DC_ell1}
\end{align}

\subsection{Result}\label{sec:main}
Recall the definitions of $\D$ and $\C$ in \eqref{eq:D_def} and \eqref{eq:C_def}. 

\subsubsection{Regime of operation}\label{sec:regime}
Our results hold in the asymptotic \emph{linear regime}, where $m,n$ and $\D$ all grow to infinity such that $m/n\rightarrow\delta\in(0,\infty)$ and $(1-\eps)m>\min_{\tau\geq 0}D(\tau)>\eps m$ for constant $\eps>0$. The assumption $m>\min_{\tau\geq 0}D(\tau)$ is motivated by \eqref{eq:m>D}.

\subsubsection{Preliminaries}
\begin{defn}[$\map$] \label{def:map}
Let $\Rc:=\{ \tau>0 | m-\D>\max\{0,\C\}\}$ and define $\map:\Rc\rightarrow(0,\infty):$
\begin{align}\label{eq:map}
\map(\tau) := \tau\frac{m-\D-\C}{\sqrt{m-\D}}.
\end{align}
\end{defn}
The next lemma shows that the inverse of $\map$ is well defined.

\begin{lem}[$\map^{-1}$,\cite{OTH}] \label{lem:map}
Assume $m>\min_{\tau>0}\D$. Then, $\Rc$ is a nonempty open interval and $\map$ is strictly increasing, continuous and bijective.  In particular, its inverse function $\map^{-1}:(0,\infty)\rightarrow\Rc$ is well defined.
\end{lem}

\subsubsection{Theorem}
Recall the assumptions of Section \ref{sec:ass}. Assume a large system setup as in Section \ref{sec:large} under the linear regime.
Theorem \ref{thm:main} characterizes the limiting behavior of the asymptotic normalized squared error of \eqref{eq:intro_GL}.
\begin{thm}\label{thm:main}
Fix any $\la>0$ in \eqref{eq:intro_GL} and let
$$\aNSE:=\lim_{\sigma\rightarrow 0}\NSE(\sigma)=\lim_{\sigma\rightarrow 0}\frac{\|\hat\x-\x_0\|_2^2}{\sigma^2}.$$
The following limit holds in probability
\begin{align}
\lim_{n\rightarrow\infty} \aNSE= \frac{D(\map^{-1}(\la))}{m-D(\map^{-1}(\la))} =:\eta(\la). \nn
\end{align}
\end{thm}

\subsection{Remarks}\label{sec:rem}

\subsubsection{The role of the parameters} Theorem \ref{thm:main} explicitly captures the role of the number of measurements $m$, the regularizer $f$, the unknown signal $\x_0$ and the regularizer parameter $\la$. The dependence on the ambient dimension $n$ is implicit through $\x_0$.

 
 \subsubsection{The mapping} The theorem maps the regularizer parameter $\la>0$ to some value $\tau\in\Rc$ through $\map^{-1}$. Note that $\Rc$ is \emph{nonempty} as long as $m>\min_{\tau} D(\tau)$ (Lemma \ref{lem:map}). Figure \ref{fig:map} illustrates the action of $\map^{-1}$ for an instance of a sparse recovery problem. 
 
 \subsubsection{Geometric nature} 
 The structure induced by $f$, the particular $\x_0$ we are trying to recover and the value of $\la$ are all summarized in a single  parameter, namely, the gaussian squared distance to the subdifferential. 
 
 
\subsubsection{Generality} In principle, Theorem \ref{thm:main} holds for any \emph{convex} regularizer $f$. Thus, it applies to any signal class that exhibits some sort
of low-dimensionality. In this sense,  it extends to the noisy case the unifying treatment of convex regularizers, which has been adopted in the analysis of noiseless compressive sensing \cite{Cha,TroppEdge}.


 \begin{figure}[!h]
\centering
  \includegraphics[width=0.9\linewidth]{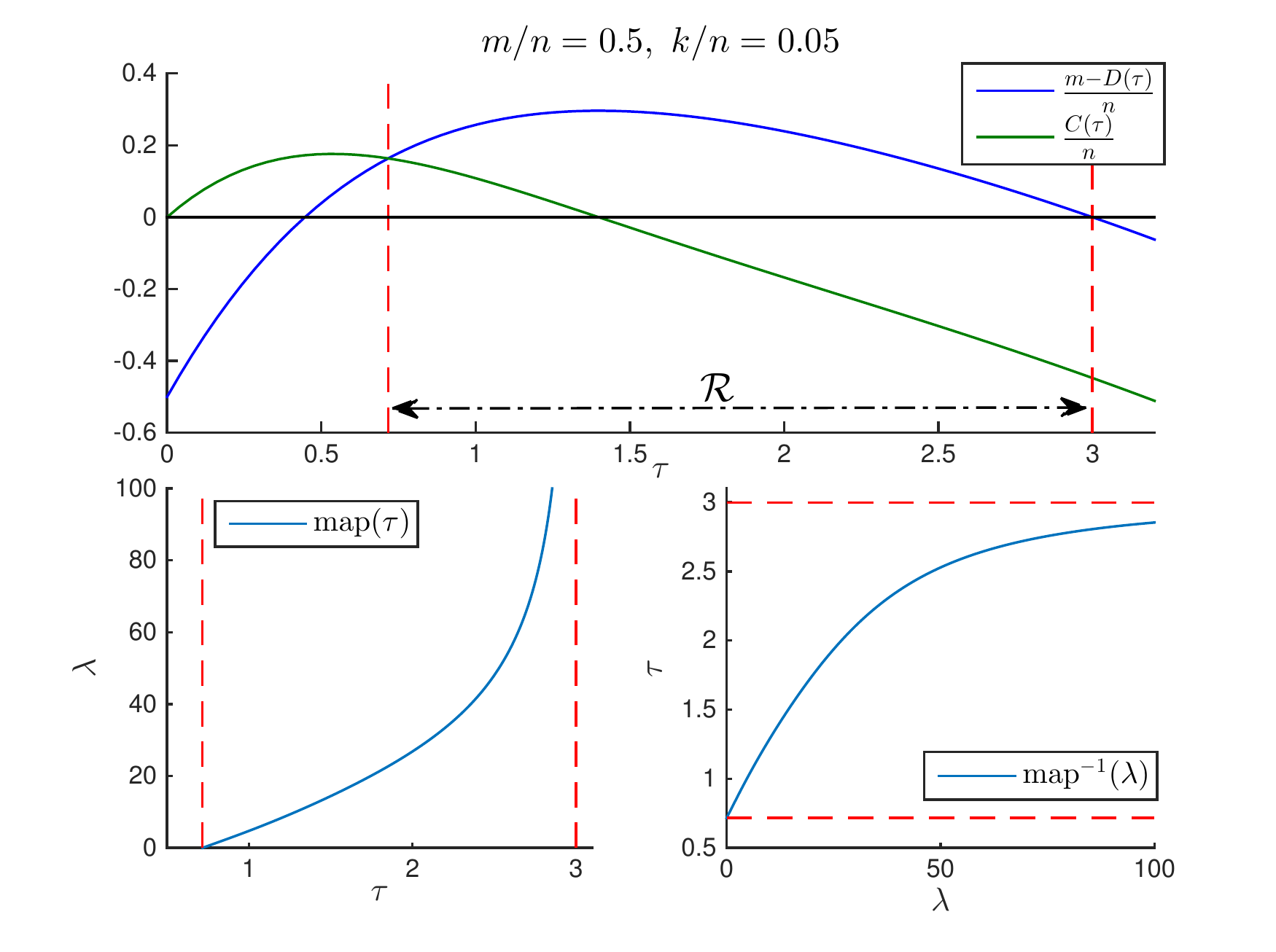}
    \label{fig:map}
  \caption{\footnotesize{Illustration of the region $\Rc$ and of the $\map$ function (Defn.~2.2) for $f=\|\cdot\|_1$ and $\x_0\in\R^{n}$ a $k$-sparse vector. $\map^{-1}$ maps the value of the regularizer $\la$ in \eqref{eq:intro_GL} to a value in $\Rc$. $D(\tau)$ and $C(\tau)$ are computed as in \eqref{eq:DC_ell1} }}.
  \label{fig:map}
\end{figure}

 \begin{figure}[!h]
\centering
  \includegraphics[width=0.85\linewidth]{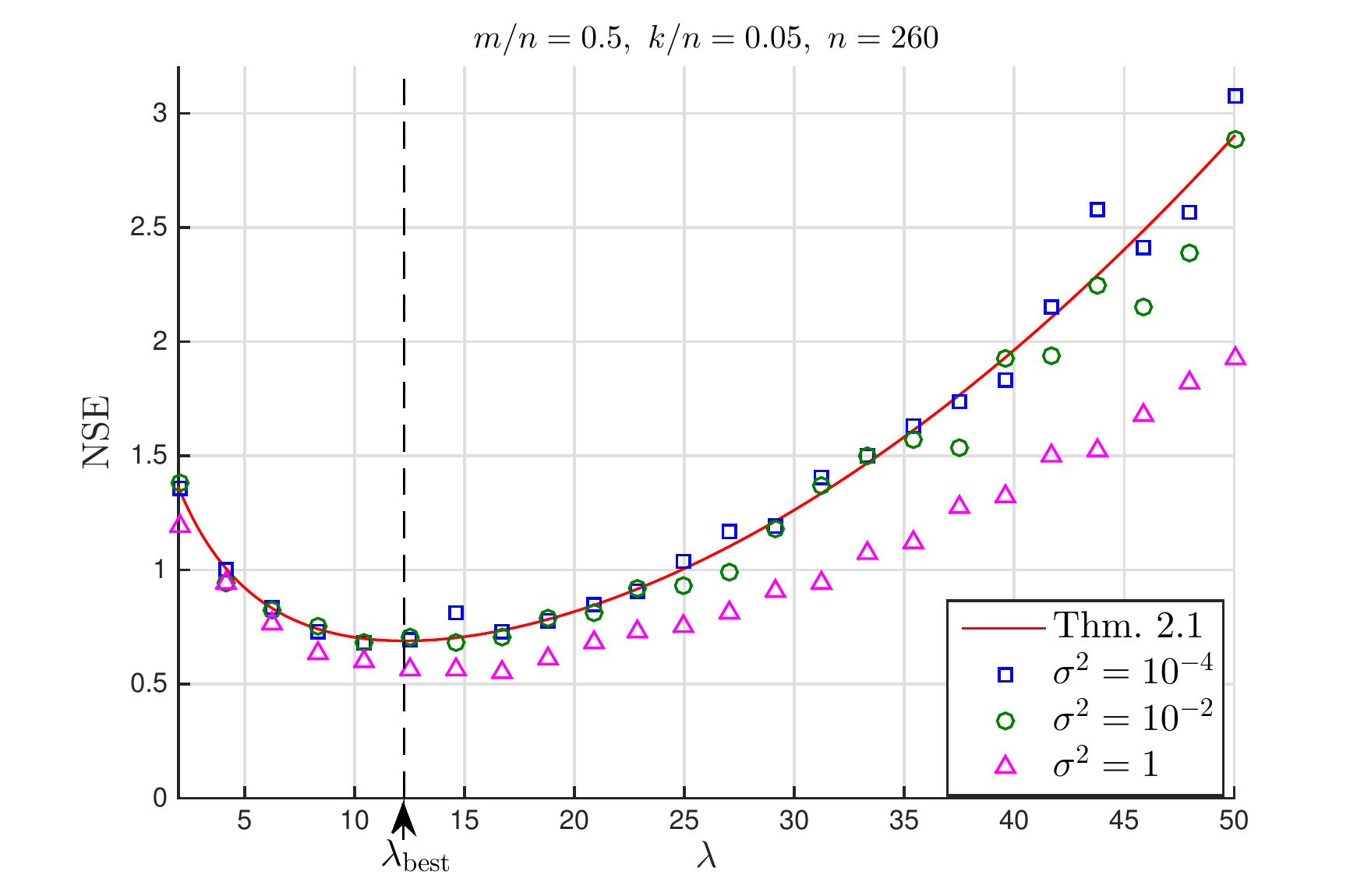}
    \label{fig:NSE}
  \caption{\footnotesize{Numerical validation of Theorem \ref{thm:main} for $f=\|\cdot\|_1$ and $\x_0\in\R^{n}$ a $k$-sparse vector. Measured values of $\NSE(\sigma)$ are averages over 50 realizations of $\A, \vb$. The theorem accurately predicts $\NSE(\sigma)$ as $\sigma\rightarrow 0$. Results support our claim that $\aNSE=\wNSE$. $\la_{best}$ is the value of the optimal regularizer as predicted by Lemma~2.2. }}
  \label{fig:NSE}
\end{figure}


\subsubsection{On the worst-case NSE} We conjecture that
\begin{align}\label{eq:conj}
\wNSE:=\sup_{\sigma>0}\NSE(\sigma)= \lim_{\sigma\rightarrow 0}\NSE(\sigma)=:\aNSE.
\end{align}
Theorem \ref{thm:main} would then imply that for any $\sigma>0$:
$$
\nlim\NSE(\sigma) \leq \eta(\la),
$$
in probability. There are several reasons that suggest this claim. First, $\wNSE=\aNSE$ has already been shown to hold for algorithms similar to \eqref{eq:intro_GL} such as: i) the constrained generalized LASSO in \eqref{eq:intro_CL}, \cite{arianSensitivity,StoLASSO,OTH}, 
ii) the proximal denoiser \cite{Don94,oymakProx}, which is essentially \eqref{eq:intro_GL} when $m=n$ and $\A=\mathbf{I}_n$. Furthermore, our conjecture is supported by  computational experiments; see Figure \ref{fig:NSE} and \cite[Sec.~13]{OTH}.

 \subsubsection{Evaluating the bound} Evaluating the bound of Theorem \ref{thm:main} for particular instances of structures and regularizers requires the ability to compute $D(\tau)$ and $C(\tau)$. It is important to note that this only requires knowledge of the particular structure of the unknown signal $\x_0$, and not the explicit unknown signal itself. For example, in sparse recovery, $D(\la)$ is the same for all $k$-sparse signals (see \eqref{eq:DC_ell1} and Fig.~\ref{fig:map}).

\subsubsection{Optimal tuning} Thm. \ref{thm:main} suggests a simple recipe for finding the optimal value $\labb$ of the regularizer parameter.
\begin{lem}\label{lem:labb}
Recall $\eta(\la)$ as defined in Theorem \ref{thm:main}.  Let $\labb:=\arg\min_{\la\geq 0}\eta(\la)$ and $\taubb:=\arg\min_{\tau\geq0}D(\tau)$. Then,
$
\labb = \taubb\sqrt{m-D(\taubb)}.
$
\end{lem}
The proof of the lemma is not involved and is omitted for brevity.
It is shown in \cite[Lem.~C.2]{TroppEdge} that $D(\tau)$ is strictly convex. Thus, $\taubb$  can be efficiently calculated as the unique solutions to a convex program. This determines $\labb$. Note that even though calculating $\labb$ does not require explicit knowledge of $\x_0$ itself, it does assume knowledge of the particular structure. For instance, in sparse recovery we need to know the sparsity level $k$ (see Fig.~\ref{fig:NSE}).


\subsubsection{Phase-transitions} Combining Theorem \ref{thm:main} with Lemma \ref{lem:labb} it holds with probability one that,
$$
\lim_{\sigma\rightarrow 0} \min_{\la>0} \frac{\|\hat\x-\x_0\|_2^2}{\sigma^2} = \frac{\min_{\tau}D(\tau)}{m-\min_{\tau}D(\tau)}.
$$
In view of the $\wNSE$ conjecture in \eqref{eq:conj}, the quantity in the left hand side can be viewed as the \emph{minimax} NSE of G-LASSO for a \emph{fixed} signal $\x_0$. While $m>\min_{\tau}D(\tau)$, we can always tune \eqref{eq:intro_GL} to guarantee robust recovery. However, as the number of measurements $m$ approaches $\min_{\tau}D(\tau)$, then, even after optimal tuning, the NSE grows to $\infty$. This phase-transition characterizing the robustness of \eqref{eq:intro_GL} is identical to \eqref{eq:m>D}, i.e. the phase-transition in noiseless compressed sensing. This observation was first formally predicted in \cite{arianSensitivity}, and, later proved in \cite{MontanariLASSO} and \cite{StoLASSO}, for $f=\|\cdot\|_1$ and $\x_0$ k-sparse. 


\subsubsection{Robustness} Theorem \ref{thm:main} reveals the following interesting feature of \eqref{eq:intro_GL}. Given sufficient number of measurements $m>\min_{\tau}D(\tau)$, the recovery is \emph{robust} for all choices of the regularizer parameter $\la>0$. In particular, this is in contrast to the $\ell_2$-LASSO in \eqref{eq:intro_GLOTH}. It was shown in \cite{OTH,TOH14} that the NSE of the later becomes unbounded if the regularizer parameter is larger than some $\la_{\max}$.


\subsubsection{Relevant literature} Most error bounds derived in the literature for \eqref{eq:intro_GL} are order-wise. The first precise results were derived in the context of sparse recovery via the AMP framework:
\cite{arianSensitivity}  develops formal expressions for the $\wNSE$ of \eqref{eq:intro_GL} under optimal tuning of the regularizer parameter $\la>0$; \cite{MontanariLASSO} explicitly characterize $\NSE(\sigma)$ for all values of $\la>0$ and all $\sigma>0$. The rest of the works that we list here use the GMT framework.  \cite{StoLASSO,OTH} precisely characterizes the $\wNSE$ of \eqref{eq:intro_CL}. \cite{OTH} computes the $\aNSE$ of \eqref{eq:intro_GLOTH}. The $\NSE(\sigma)$ of \eqref{eq:intro_GLOTH} with $\ell_1$-regularization but arbitrary $\sigma>0$ has been characterized by the authors in \cite{ICASSP}.  Theorem \ref{thm:main} characterizes the $\aNSE$ of the generalized $\ell_2^2$-LASSO.

\section{Proof Outline}\label{sec:proof}
We outline the main steps of the proof here. Most of the technical details are deferred to the Appendix.
Before everything, we re-write \eqref{eq:intro_GL} by changing the decision variable to be the error vector $\ww=\x-\x_0$:
\begin{align}\label{eq:GL}
\hat\ww := \min_{\ww}\frac{1}{2}\|\A\ww - \sigma\vb \|_2^2 + \frac{\la}{\sigma}f(\x_0 + \ww).
\end{align}
Theorem \ref{thm:main} states a precise expression for the limiting behavior $\lim_{\sigma\rightarrow 0}\|\hat\ww\|^2/\sigma^2$.  Throughout the analysis, we fix any $\la>0$. Also, we simply write $\|\cdot\|$ instead of $\|\cdot\|_2$.

\subsection{First-order Approximation}
We start with a useful approximation to \eqref{eq:GL}.
 The idea is that in the regime of interest we expect $\hat\ww$ to scale linearly with $\sigma$. Thus, in the limit $\sigma\rightarrow 0$, $\|\hat\ww\|$ is sufficiently small such that $f(\x_0+\ww)\approx f(\x_0) + \sup_{\s\in\paf}\s^T\ww$. Note that this always holds  with a ``$\geq$" sign due to convexity. What we show in the Appendix is essentially that introducing this approximation in \eqref{eq:GL} does note alter $\|\hat\ww\|$ in the limit $\sigma\rightarrow 0$.

\subsection{Gaussian min-max Theorem}\label{sec:GMT}
We get a handle on \eqref{eq:GL} and its optimal value via analyzing a different and simpler optimization problem. The machinery that allows this relies on  Gordon's Gaussian min-max theorem (GMT) \cite[Lem.~3.1]{gorLem}. In fact, we require a stronger version of the GMT that can be obtained when accompanied with additional \emph{convexity} assumptions that are not present in its original formulation. The fundamental idea is attributed to Stojnic \cite{StoLASSO}. \cite{tight} builds upon this and derives a concrete and somewhat extended statement of the result in \cite[Thm.~II.1]{tight}. Please refer to the discussion in \cite{tight} for further details on the GMT, the role of convexity, and, the differences between \cite[Lem.~3.1]{gorLem}, \cite{StoLASSO} and \cite[Thm.~II.1]{tight}. We summarize the result of \cite[Thm.~II.1]{tight} in the next few lines.
Let $\G\in\R^{m\times n},\g\in\R^m,\h\in\R^n$ have entries i.i.d. Gaussian; $\Sc_\ab\subset\R^n$,$\Sc_\bb\subset\R^m$ be convex compact sets, and $\psi:\Sc_\ab\times\Sc_\bb\rightarrow\R$ be convex-concave and continuous. Further consider the following two min-max problems:
\begin{align}\label{eq:Phi}
\Phi(\A) := \min_{\ab\in\Sc_\ab} \max_{\bb\in\Sc_\bb} \bb^T\G\ab + \psi(\ab,\bb),
\end{align}
\begin{align}{\label{eq:phi}
\phi(\g,\h) := \min_{\ab\in\Sc_\ab} \max_{\bb\in\Sc_\bb} \|\ab\|\g^T\bb - \|\bb\|\h^T\ab + \psi(\ab,\bb).}
\end{align}
Then, for any $\mu\in\R, t>0$:
$$
\Pro\left(|\Phi(\A) - \mu|>t\right) \leq 2 \Pro\left(|\phi(\g,\h) - \mu|>t\right).
$$
Thus, if the optimal cost $\phi(\g,\h)$ of \eqref{eq:phi} concentrates to some value $\mu$, the same is true for $\Phi(\A)$. This suggests analyzing \eqref{eq:phi} instead of \eqref{eq:Phi}, and indirectly yield conclusions for the latter. The premise is that the optimization in \eqref{eq:phi} is easier to analyze; we often refer to it as ``Gordon's optimization" following \cite{tight}. Assuming a setup in which the problem dimensions $m,n$ grow to infinity it is shown in \cite{tight} that if $\phi(\g,\h)$ converges in probability to deterministic value $d_*$, then, so is $\Phi(\G)$. What is more, if $\|\ab_*(\g,\h)\|$ converges to say $\alpha_*$, and some appropriate strong convexity assumption on the objective function of \eqref{eq:phi} is satisfied, then $\|\ab_*(\G)\|$ also converges to $\alpha_*$.    Here, we have denoted $\ab_*(\g,\h)$, $\ab_*(\G)$ for the minimizers in \eqref{eq:phi} and \eqref{eq:Phi}, respectively; refer to \cite{tight} and Lemma \ref{lem:tight} for the exact statements. As might be already suspected, this latter property is of interest to our problem. In what follows, we bring \eqref{eq:intro_GL} in the format of \eqref{eq:Phi}, derive the corresponding ``Gordon's optimization" problem and analyze the minimizer of that one instead.

\subsection{Gordon's Optimization}
We use the fact that 
$(1/2)\|\ab\|^2 = \max_{\bb} \bb^T\ab - (1/2)\|\bb\|^2$ 
to equivalently express $\hat\ww$ as the solution to (also, recall the first-order approximation)
\begin{align}\nn
\min_{\ww}\max_{\bb}\bb^T\A\ww - \sigma\bb^T\vb -({1}/{2})\|\bb\|^2 + \la\max_{\s\in\paf}\s^T\ww.
\end{align}
Identify { $\psi(\ww,\bb)=- \sigma\bb^T\vb -\frac{1}{2}\|\bb\|^2 + \la\max_{\s}\s^T\ww$}, which is convex-concave and continuous, to see that the above is in the desired format \eqref{eq:Phi}. The only caveat is that the constraint sets on $\ww$ and $\bb$ appear unbounded. This is appropriately treated in the Appendix and we do not elaborate any further here. The corresponding ``Gordon's optimization" problem writes:
\begin{align*}
\min_{\ww}\max_{{{ \s, \bb}  }}\|\ww\|\g^T\bb-\|\bb\|\h^T\ww - \sigma\bb^T\vb -\frac{1}{2}\|\bb\|^2 + \la\s^T\ww,
\end{align*}
where the variable $\s$ is constrained in $\paf$, but we omit to shorten notation. In the next lines, we show how to simplify this optimization to a scalar problem. 
Recall that $\g,\vb$ both have entries i.i.d. $\Nn(0,1)$ and are independent of each other; thus $\|\ww\|\g+\sigma\vb$ has entries $\Nn(0,\sqrt{\|\ww\|^2+\sigma^2})$. Also, note that the maximization over the direction of $\bb$ is easy to perform, since $\max_{\|\bb\|=\beta}\g^T\bb=\beta\|\g\|$, $\beta\geq0$. With these, and some abuse of notation so that $\g$ continues being i.i.d standard normal gaussian, we may rewrite the above  as:
\begin{align}\nn
\min_{\ww}\max_{\s, \beta\geq 0} \sqrt{\|\ww\|^2+\sigma^2}\|\g\|\beta-(\beta\h-\la\s)^T\ww -{\beta^2/2}.
\end{align}
Observe that the objective function is now convex in $\ww$ and concave in $\beta,\s$. Thus, modulo compactness of the constraint sets (see Appendix for details), we can flip the order of minimization and maximization \cite[Cor.~37.3.2]{Roc70}, and write:
\begin{align}\nn
\max_{\s, \beta\geq 0}\min_{\ww} \sqrt{\|\ww\|^2+\sigma^2}\|\g\|\beta-(\beta\h-\la\s)^T\ww -{\beta^2/2}.
\end{align}
But, now, it is easy to perform the minimization over the \emph{direction} of $\ww$. Doing this, and letting $\alpha$ represent $\|\ww\|$:
\begin{align}\nn
\max_{\s, \beta\geq 0}\min_{\alpha\geq 0} \sqrt{\alpha^2+\sigma^2}\|\g\|\beta-\alpha\beta\|\h-\frac{\la}{\beta}\s\| -{\beta^2/2}.
\end{align}
We are almost done with the simplifications. One last step amounts to flipping the order of min-max once more (the objective is appropriately concave-convex) and performing the maximization over $\s$, which results in the appearance of the distance term below:
\begin{align}{\label{eq:GO}
\min_{\alpha\geq 0}\max_{\beta\geq 0} \sqrt{\alpha^2+\sigma^2}\|\g\|\beta-\alpha\dt(\h,\frac{\la}{\beta}\paf) -{\beta^2/2}.}
\end{align}

%
\subsection{Analysis of Gordon's Optimization}
In \eqref{eq:GO}, the variable $\alpha$ plays the role of $\|\ww\|$. Thus, from the discussion in Section \ref{sec:GMT}, if we find the value to which the optimal $\alpha_*(\g,\h)$ in \eqref{eq:GO} converges, then, we may conclude that the desired quantity $\|\hat\ww(\A,\vb)\|$ also converges to the same value. This will establish Theorem \ref{thm:main}. 
Assume the asymptotic regime that holds for Theorem \ref{thm:main}
We only highlight the main ideas here and defer most of the details to the Appendix.

Both functions $\|\g\|$ and $\dt(\h,(\la/\beta)\paf)$ in \eqref{eq:GO} are 1-Lipschitz in their arguments. Then, the classical gaussian concentration of Lipschitz functions implies that they concentrate around $\sqrt{m}$ and $\sqrt{D(\la/\beta)}$, respectively (e.g. \cite[Lem.~B.2]{OTH}).  We use this in the appendix to prove that \eqref{eq:GO} converges in probability (after proper normalization) to
\begin{align}\label{eq:d}
\min_{\alpha\geq 0}\max_{\beta\geq 0} \sqrt{\alpha^2+\sigma^2}\sqrt{m}\beta-\alpha\beta\sqrt{D(\la/\beta)}-{\beta^2/2}.
\end{align}
Moreover, the minimizer of \eqref{eq:GO} converges to the minimizer $\alpha_*$ of the deterministic minimization in \eqref{eq:d}.  To compute $\alpha_*$, we use duality (the objective is (strictly) convex in $\alpha$ and concave in $\beta$). First, fix $\beta$, differentiate the objective in \eqref{eq:d} w.r.t. $\alpha$ and equate to $0$ to find that is minimized at 
$$
\alpha_*(\beta) = \sqrt{D(\la/\beta)}/\sqrt{m-D(\la/\beta)}.
$$
Substituting this value back in \eqref{eq:d} and differentiating now with respect to $\beta$, yields the optimal $\beta_*=\la/\map^{-1}(\la)$. Note that $\alpha_*(\beta_*)$ agrees with the expression of the theorem to conclude.



\bibliography{\myreferences}

\appendix

In Section \ref{sec:proof} we outlined the proof of Theorem \ref{thm:main}. Here, we provide a complete proof of the theorem.

\subsection{Preliminaries}\label{sec:pre}
We rewrite \eqref{eq:intro_GL} in a more convenient format for the purposes of the analysis. In particular, (i) substitute $\y=\A\x_0+\sigma\vb$, (iii) subtract from the objective the constant term $\la f(\x_0)$,  (ii) change the decision variable to the quantity of
interest, i.e. the normalized error vector $\w:=(1/\sigma)({\x-\x_0})$, (iv) rescale by a factor of $\sigma^2$. Then,
\begin{align}\label{eq:GL_app_app}
\wh := \min_{\w}\frac{1}{2}\|\A\w - \vb \|_2^2 + \frac{\la}{\sigma}( f(\x_0 + \sigma \w) - f(\x_0)).
\end{align}
We will derive a precise expression for the limiting (as $n\rightarrow\infty$) behavior of $\lim_{\sigma\rightarrow 0}\|\wh\|_2$. Note that after the normalization of $\x-\x_0$ with $\sigma$ , it is not guaranteed that the optimal minimizer in \eqref{eq:GL_app_app} is bounded (think of $\sigma\rightarrow 0$). However, we will prove that in the regime of Theorem \ref{thm:main} this is indeed the case. Many of the arguments that we use in the analysis require boundedness of the constraint sets. 
To tackle this, we assume that $\wh$ is bounded by some large constant $K>0$ (with probability one over $\A,\vb$), the value of which to be chosen at the end of the analysis. Recall that at that point we will have a precise characterization of the limiting behavior of $\|\wh\|_2$, say $\alpha_*$. If $\alpha_*$ turns out to be independent on the value of $K$ which we started with, then we will assume that this starting value was strictly larger than $\alpha_*$.  Thus, in what follows, we let $K,\La,M$ denote such (arbitrarily) large, but finite, positive quantities. For $K$, which is reserved as an upper bound on $\|\w\|$ we assume that is constant in the sense that it does not scale with $n$. This will be required when we apply \cite[Thm.~II.2]{tight} in Section \ref{sec:applyThm}. On the other hand, $\La,\M$ are in general allowed to depend on $n$.
Also, we fix $\la>0$ and write $\|\cdot\|$ instead of $\|\cdot\|_2$.

\subsection{Gordon's Optimization for arbitrary $\sigma$}
We use the fact that 
\begin{align}\label{eq:trick}
(1/2)\|\ab\|^2 = \max_{\bb} \bb^T\ab - (1/2)\|\bb\|^2,
\end{align}
to equivalently express $\wh$ as the solution to 
\begin{align}
\min_{\|\w\|\leq K}\max_{\|\bb\|\leq\Lambda} \bb^T\A\w & - \bb^T\vb - \frac{1}{2}\|\bb\|^2\nn
 \\ &+ \frac{\la}{\sigma}( f(\x_0+\sigma\w) - f(\x_0)).\nn
\end{align}
In view of \eqref{eq:trick} and the boundedness of $\w$, the set of optima of $\bb$ is also bounded by some $0<\Lambda:=\Lambda(K)<\infty$.
This brings \eqref{eq:GL_app_app} in the desired format \eqref{eq:Phi}. Œ
 Then, \eqref{eq:phi} writes
\begin{align}
\wt(\g,\h) &:= \min_{\|\w\|\leq K}\max_{\|\bb\|\leq\Lambda} \sqrt{\|\w\|^2+1}\g^T\bb  \nn
- \|\bb\|\h^T\w \\&\qquad- \|\bb\|^2/2 + ({\la}/{\sigma})( f(\x_0+\sigma\w) - f(\x_0) )\nn,
\end{align}
The maximization over the direction of $\bb$ is easy to perform; note that $\max_{\|\bb\|=\beta}\g^T\bb=\beta\|\g\|_2, \beta\geq 0$. Also,
$f$ is continuous and convex, thus, we can express it in terms of its convex conjugate $f^*(\ub) = \sup_{\x}\ub^T\x - f(\x)$. In particular, applying \cite[Thm.12.2]{Roc70} we have $f(\x_0+\sigma\w)=\sup_{\ub}\x_0^T\ub + \sigma\ub^T\w-f^*(\ub)$. The supremum here is achieved at $\ub_*\in\partial f(\x_0+\sigma\w)$ \cite[Thm.~23.5]{Roc70}. Also, from \cite[Prop.~4.2.3]{Bertsekas}, $\cup_{\|\w\|\leq K}\partial f(\x_0+\sigma\w)$ is bounded. Thus, the set of maximizers $\ub_*$ is bounded and for some $0< M:=M(K) <\infty$, $\wt$ is given as the solution to
\begin{align}\label{eq:GO2}
&\phi(\sigma;\g,\h):=\min_{\|\w\|\leq K}\max_{\substack{0\leq\beta\leq\La \\ 0\leq \|\ub\| \leq M}  }
\sqrt{\|\w\|^2+1}\|\g\|\beta - {\beta^2}/{2} \nn \\ &~~~ -(\beta\h-\la\ub)^T\w  + \frac{\la}{\sigma}(\ub^T\x_0 - f^*(\ub)-f(\x_0)).
\end{align}

 \subsection{Gordon's Optimization in the limit $\sigma\rightarrow 0$}

\cite[Thm.~II.2]{tight} relates $\|\wt\|$ to $\|\wh\|$, under appropriate assumptions. Also, recall that  we wish to characterize $\lim_{\sigma\rightarrow 0}\|\wh\|$. Thus, in view of \eqref{eq:GO2} we wish to analyze the problem $$\phiz:=\phiz(\g,\h):=\lim_{\sigma\rightarrow 0}\phi(\sigma;\g,\h).$$ In \eqref{eq:GO2}, from Fenchel's inequality:
\begin{align}\label{eq:conj}
\ub^T\x_0 - f^*(\ub)-f(\x_0)\leq 0.
\end{align}
With this observation, we prove in the next lemma that $\phi(\sigma;\g,\h)$ is non-decreasing in $\sigma$; see Section \ref{sec:lemproof} for the proof.
\begin{lem}\label{lem:decreasing}
Fix $\g,\h$ and consider $\phi(\cdot;\g,\h):(0,\infty)\rightarrow\R$ as defined in \eqref{eq:GO2}. $\phi(\sigma;\g,\h)$ is non-decreasing in $\sigma$.
\end{lem}
%
In particular, when viewed as a function of $\kappa:=\la/\sigma$, $\phi(\cdot;\g,\h)$ is non-increasing. Thus,
 \begin{align}\label{eq:liminf}
 \phiz = \lims\phi(\sigma)= \lim_{\kappa\rightarrow\infty}\phi(\kappa)=\inf_{\kappa\geq 0}\phi(\kappa),
 \end{align}
Next, we argue that we can flip the order of min-max; we will apply  \cite[Cor.~37.3.2]{Roc70}.
 The objective function in \eqref{eq:GO2} is continuous,  convex in both $\w,\kappa$, and, concave both in $\bb,\ub$. The constraint sets are all convex and one of them is bounded. With this and \eqref{eq:liminf}, we get
%
\begin{align}\nn
&\phi_0(\g,\h)=\max_{\substack{0\leq \beta\leq \Lambda \\ 0\leq\|\ub\|\leq M}}\min_{\|\w\|\leq K }\inf_{\kappa\geq 0}  
\sqrt{\|\w\|^2+1}\|\g\|\beta- \frac{\beta^2}{2} \\ &~~~~ -(\beta\h-\la\ub)^T\w  + \kappa(\ub^T\x_0 - f^*(\ub)-f(\x_0)).\label{eq:sigma}
\end{align}
Recall \eqref{eq:conj} and the fact that equality is achieved iff $\ub\in\paf$ (e.g. \cite[Thm.~23.5]{Roc70}). Then, $\phi_0(\g,\h)$ is given by
\begin{align*}
\max_{\substack{0\leq\beta\leq \Lambda \\ \ub\in\paf}}\min_{\|\w\|\leq K}  
\sqrt{\|\w\|^2+1}\|\g\|\beta \nn - (\beta\h-\la\ub)^T\w - \frac{\beta^2}{2},
\end{align*}
where we have assumed $\infty>M>\max_{\s\in\paf}\|\s\|$.
We can simplify this one step further by performing the minimization over the direction of $\w$. In the problem below note that $\alpha$ plays the role of $\|\w\|$. Thus,
$\phi_0(\g,\h) = \max_{\substack{0\leq\beta\leq \Lambda \\ \ub\in\paf}}\min_{0\leq\alpha\leq K}  
\sqrt{\alpha^2+1}\|\g\|\beta \nn - \alpha\|\beta\h-\la\ub\| - \frac{\beta^2}{2}.$
The objective function  above is continuous, convex in $\alpha$ and concave in $\beta,\ub$. Also the constraint sets are convex and bounded. Thus,  \cite[Cor.~37.3.2]{Roc70}, we can flip the order of max-min. 
Also, for $\beta>0$, $\min_{\ub\in\paf}\|\beta\h-\la\ub\| = \beta\dt(\h,(\la/\beta)\paf)$. With these, normalizing with $1/m$ and  appropriately rescaling $\beta$:
\begin{align}\label{eq:goropt}
&\tilde\phi_0(\g,\h) = \min_{0\leq \alpha\leq K}\max_{0\leq\beta\leq \Lambda}  \Lc(\alpha,\beta;\g,\h):= \nn\\
&\hspace{-20pt}\quad\sqrt{\alpha^2+1}\frac{\|\g\|}{\sqrt{m}}\beta - \alpha\beta\frac{1}{\sqrt{m}}\dt(\h,\frac{\la}{\beta\sqrt{m}}\paf) - \frac{\beta^2}{2}.
\end{align}
Normalization here is convenient for the purposes of applying statement (iii) of \cite[Thm.~II.2]{tight}, which follows.

\subsection{Applying \cite[Thm.~II.2]{tight}}\label{sec:applyThm}

\eqref{eq:goropt} describes ``Gordon's optimization" (modulo normalization with $m$) corresponding to \eqref{eq:GL_app_app} in the limit of $\sigma\rightarrow 0$. Also, the variable $\alpha$ in \eqref{eq:goropt} plays the role of $\|\w\|$. The idea is that problem \eqref{eq:goropt} behaves in the large system limit just like \eqref{eq:GL_app_app}. This is formalized in the lemma below, which is a direct corollary of \cite[Thm.~II.2]{tight} applied to our setup and combined with our analysis thus far. 
For the statement of the lemma recall that we are operating in the large-system limit in which problem dimensions $n,m,\min_{\tau\geq 0}D(\tau)$ grow \emph{linearly} to infinity (cf. Section \ref{sec:regime}). Also, we use standard notation $X_n\rP c$, to denote convergence in probability of $\{X_n\}_{n=1}^\infty$ to $c\in\R$ as $n\rightarrow\infty$.

\begin{lem}\label{lem:tight}
Let $\wh$ and $\cost_\sigma$ be a minimizer and the optimal cost of \eqref{eq:GL_app_app}, respectively. Also, recall \eqref{eq:goropt}. Suppose $d:[0,\infty)\rightarrow\R$ is such that $\max_{0\leq\beta\leq\La}\Lc(\alpha,\beta;\g,\h)\rP d(\alpha)$ for all $\alpha\in[0,K]$, and $d(\alpha)\geq d(\alpha_*)+\zeta(\alpha-\alpha_*)^2,\forall\alpha\in[0,K]$, for some $\alpha_*\in[0,K]$ and $\zeta>0$. Then,
\begin{align}
\lim_{\sigma\rightarrow 0}\|\wh\|_2\rP \alpha_* \quad \text{ and } \quad \lim_{\sigma\rightarrow 0}\frac{\cost_\sigma}{\sigma^2 m}\rP d(\alpha_*).
\end{align}
\end{lem}

Recall $\wh= (\hat\x - \x_0)/\sigma$, thus, $\|\wh\|^2=\NSE(\sigma)$. 
In what follows, we construct deterministic function $d$ that satisfies the requirements of Lemma \ref{lem:tight} and prove that $\alpha_*$ satisfies the formula of Theorem \ref{thm:main}. This will complete the proof. To suppress notation, define $\labe:=\la/(\beta\sqrt{m})$.

Both functions $\|\g\|$ and $\dt(\h,\labe\paf)$ are 1-Lipschitz in their arguments. Then, the classical gaussian concentration of Lipschitz functions implies that they concentrate around $\sqrt{m}$ and $\sqrt{D\labe}$, respectively 

From standard concentration results on Lipschitz functions of gaussian r.v.s.  (e.g. \cite[Lem.~B.2]{OTH}), we have for all $\alpha,\beta$:
\begin{align}\label{eq:ptw}
\Lc(\alpha,\beta;\g,\h)\rP \dd(\alpha,\beta):= \sqrt{\alpha^2+1}\beta - \alpha\beta\sqrt\frac{D\labe}{m} - \frac{\beta^2}{2}.
\end{align}
As we have seen $\Lc(\alpha,\beta;\g,\h)$ is convex in $\alpha$ and concave in $\beta$. Since taking limits preserves convexity, the same is true for $\dd(\alpha,\beta)$. Next, define $$d(\alpha):=\max_{0\leq\beta\leq \La}\dd(\alpha,\beta).$$ 
We claim that this satisfies the prerequisites of Lemma \ref{lem:tight}. 

First, we show the convergence part. It suffices to prove that for each $\alpha$ the convergence in \eqref{eq:ptw} holds uniformly over all $\beta\in[0,\La]$. Concavity of $\Lc(\alpha,\beta)$ in its second argument is critical. In particular, the claim follows from \cite[Cor.~II.1]{AG1982}: ``point-wise convergence in probability of concave functions implies uniform convergence in compact spaces".

Next, we compute appropriate $\alpha_*$. 
Consider 
\begin{align}\label{eq:dab}
(\alpha_*,\beta_*) := \arg\min_{0\leq\alpha\leq K}\max_{0\leq\beta\leq \Lambda}\dd(\alpha,\beta)
\end{align}
 We compute those in the next lemma; see Section \ref{sec:lemproof} for a proof.

\begin{lem} \label{lem:det}
Consider the optimization in \eqref{eq:dab}. Let
$$
 \alpha_* =\sqrt{\frac{D(\map^{-1}(\la))}{m-D(\map^{-1}(\la))}}  \text{ and } \beta_* = \frac{\la}{\map^{-1}(\la)\sqrt{m}} .
$$
Then, there exist $K,\La$ satisfying $0<\alpha_*<K<\infty$ and $0<\beta_*<\La<\infty$, such that $(\alpha_*,\beta_*)$ are optimal in \eqref{eq:dab}. 
%
\end{lem}

Let $K,\La$ in \eqref{eq:dab} be as in Lemma \ref{lem:det}.
It remains to prove $d(\alpha)>d(\alpha_*)+\zeta(\alpha-\alpha_*)^2$ for some $\zeta>0$.
Fix $\alpha\in[0,K]$. Clearly,
$$d(\alpha) = \max_{\beta} \tilde d(\alpha,\beta) \geq \tilde d(\alpha,\beta_*).$$
We will now use the fact that for fixed $\beta\in(0,\La]$, the function $\tilde d(\alpha,\beta)$ is strongly convex in $0\leq \alpha\leq K$. Indeed,
$$
{\partial^2{\tilde d}}/{\partial\alpha^2} = \beta/{(\alpha^2+1)^{3/2}} \geq \beta/{(K^2+1)^{3/2}}.
$$
Recall $\beta_*>0$ and let $\zeta = \beta_*/ {(K^2+1)^{3/2}}>0$. Then,
\begin{align}\nn
d(\alpha) \geq \tilde d(\alpha,\beta_*) \geq \tilde d(\alpha_*,\beta_*) + \zeta(\alpha-\alpha_*)^2.
\end{align}

\subsection{Proofs of Auxiliary Results}\label{sec:lemproof}

\subsubsection{Lemma \ref{lem:decreasing}}
Denote $\Lc(\sigma,\w,\beta,\ub)$ the objective function in \eqref{eq:GO2} and consider $0<\sigma_1<\sigma_2<\infty$. Let $\w^{(2)}$ be an optimal solution to the min-max problem in \eqref{eq:GO2} for $\sigma_2$. Then, let $(\beta^{(1)},\ub^{(1)}) = \arg\max_{\beta,\ub}\Lc(\sigma_1,\w^{(2)},\beta,\ub)$. Clearly,
$$
\phi(\sigma_1)\leq \Lc(\sigma_1,\w^{(2)},\beta^{(1)},\ub^{(1)}).
$$
Using $\la/\sigma_1>\la/\sigma_2$ and  \eqref{eq:conj},
$$
 \Lc(\sigma_1,\w^{(2)},\beta^{(1)},\ub^{(1)}) \leq  \Lc(\sigma_2,\w^{(2)},\beta^{(1)},\ub^{(1)})
$$
But,
$$
\Lc(\sigma_2,\w^{(2)},\beta^{(1)},\ub^{(1)}) \leq \phi(\sigma_2).
$$
Combine the above chain of inequalities to  conclude.

\subsubsection{Lemma \ref{lem:det}}
%
Let $\beta_*, \alpha_*$ be as in the statement of the lemma. Also $\tau_*:=\map^{-1}(\la)>0$ (cf. Definition \ref{def:map}). Notice that $\alpha_*,\beta_*>0$ and set 
$$
0<K=2\alpha_*<\infty \quad \text{ and } \quad 0<\La=2\beta_*<\infty.
$$
As we have seen $\tilde d(\alpha,\beta)$ is convex-concave. Also, the constraint sets are convex and compact, hence,
$$
(\alpha_*,\beta_*) = \max_{0\leq\beta\leq \Lambda}\min_{0\leq\alpha\leq K}\dd(\alpha,\beta).
$$
We have $0<\beta_*<\La$ and $D(\la_{\beta_*})<m$ and $D\labe$ continuous in $\beta$ (cf. \cite[Lem.~8.1]{OTH}). Thus, there exists open neighborhood $\Nn_1\subset[0,\La]$ such that $D\labe<m, \forall\beta\in\Nn_1$. Fix any such $\beta\in\Nn_1$ and let 
\begin{align}\label{eq:convex}
d_*(\beta): = \min_{0\leq\alpha\leq K} \tilde d(\alpha,\beta).
\end{align}
 Differentiating with respect to $\alpha$, we find
$$
\frac{\partial{\tilde d}}{{\partial\alpha}} = \beta\frac{\alpha}{\sqrt{\alpha^2+1}}- \beta\sqrt\frac{{D\labe}}{{m}}.
$$
It can be checked that $\alpha_*(\beta) = \sqrt{\frac{D\labe}{m-D\labe}}$ is the unique solution to the equation ${\partial{\tilde d}}/{{\partial\alpha}}=0$. In particular, $\alpha_*=\alpha_*(\beta_*)$ and is feasible, i.e. $\alpha_*(\beta_*)\in(0,K)$. From continuity of $D(\cdot)$ (cf. \cite[Lem.~8.1]{OTH}), we have  $\alpha_*(\beta)$ be a continuous function of $\beta$. Hence, there exists open neighborhood of $\beta_*$, say $\Nn_2\subset\Nn_1$, such that $\alpha_*(\beta)\in(0,K)$ for all $\beta\in\Nn_2$. For any such $\beta\in\Nn_2$, $\alpha_*(\beta)$ satisfies first-order optimality conditions of the convex minimization in \eqref{eq:convex}, thus, is optimal:
$$
d_*(\beta)  =  -{\beta^2}/{2} + \frac{\beta}{\sqrt{m}}\sqrt{m-D\labe}, ~~\forall\beta\in\Nn_2.
$$
Differentiating this with respect to $\beta$ finds:
$$
\frac{\partial{ d_*}}{\partial\beta} = -\beta + \frac{1}{\sqrt{m}}\frac{m-D\labe-C\labe}{\sqrt{m-D\labe}},
$$
where we have used $\partial D(\tau)/\partial\tau = -(2/\tau) \C,\tau>0$ (cf.\cite[Lem.~C.2]{TroppEdge}). Note that the second summand above is equal to $\frac{\beta}{\la}\map(\la/(\beta\sqrt{m}))$. With this, it is easy to verify that $\beta_*$ is such that $\partial{d_*}/\partial\beta=0$. From \eqref{eq:convex}, $d_*(\beta)$ is concave as the point-wise minimum of concave functions. Thus, first-order optimality conditions satisfied by $\beta_*$ are sufficient, which completes the proof.

\end{document}